\documentclass[11pt,a4paper,Color]{article}
\usepackage{amssymb}
\usepackage{amsmath}
\usepackage{graphicx, color}
\usepackage{indentfirst}
\newcommand{\mc}[1]{\mathcal{#1}}

\input xypic.sty

\def\qed{\hfill {\large ${\sqcup\!\!\!\!\sqcap}$}}

\newenvironment{demo}{{\bf Proof }}
{\qed \\}

\newcommand{\re}{\mathbb R}

\newcommand{\flecha}{\longrightarrow}
\newcommand{\<}{\left<}
\renewcommand{\(}{\left(}
\newcommand{\lb}{\label}
\newcommand{\nn}{\nonumber}
\newcommand{\fracc}{\displaystyle\frac}

\renewcommand{\>}{\right>}
\renewcommand{\)}{\right)}

\def\bal{\begin{align}}
\def\eal{\end{align}}

\newcommand{\bde}{\begin{defi}}
\newcommand{\ede}{\end{defi}}

\numberwithin{equation}{section}
\def\be{\begin{equation}}
\def\ee{\end{equation}}

\def\oH{{\overline H}}
\def\og{{\overline g}}
\def\oR{{\overline R}}
\def\oRic{{\overline Ric}}
\def\oM{{\overline M}}
\def\oN{{\overline \nabla}}

\def\t{{\frak{t}}}
\def\a{\alpha}

\def\tr{{\rm tr}}

\def\grad{{\rm grad}}

\def\nablaa{\overline{\nabla}}

\def\vle{{\rm vol}}
\def\parcial#1#2{\frac{\partial #1}{\partial#2}}

\def\flecha{\longrightarrow}

\def\ds{\displaystyle}

\newtheorem{defi}{Definition}
\newtheorem{teor}{Theorem}
\newtheorem{prop}[teor]{Proposition}

\newtheorem{lema}[teor]{Lemma}

\newtheorem{nota}{Remark}
\newtheorem{coro}[teor]{Corollary}

\numberwithin{lemap}{teor}
\numberwithin{corop}{teor}
\numberwithin{ejer}{subsection}
\numberwithin{ejemplo}{subsection}
\begin{document}

\title{Volume\,-preserving mean curvature flow of revolution hypersurfaces in a Rotationally Symmetric Space}

\pagestyle{myheadings}\markboth{E. Cabezas-Rivas and V. Miquel}{Volume-preserving M.C.F. in Rotationally Symmetric Spaces}%\markleft{}

\title{Volume\,-preserving mean curvature flow of revolution hypersurfaces in a Rotationally Symmetric Space}% \footnote{Mathematics Subject
%Classification(2000) 53C42, 52C21}}

\author{ Esther Cabezas-Rivas and
Vicente Miquel}

\date{}

\maketitle

\begin{abstract}
In an ambient space with  rotational symmetry around an axis (which include the Hyperbolic and Euclidean spaces), we study the evolution under the volume-preserving mean curvature flow of a revolution hypersurface $M$ generated by a graph over the axis of revolution and with boundary in two totally geodesic hypersurfaces (tgh for short). Requiring that, for each time $t\geq 0$, the evolving hypersurface $M_t$ meets such tgh ortogonally,  we prove that: a) the flow 
exists while $M_t$ does not touch the axis of rotation; b) throughout the time interval of existence, b1) the generating curve of $M_t$ remains a graph, and b2) the averaged mean curvature is double side bounded by positive constants; c) the singularity set (if non-empty) is finite and lies on the axis; d) under a suitable hypothesis relating the enclosed volume to the $n$-volume of $M$, we achieve long time existence and convergence to a revolution hypersurface of constant mean curvature. 
\end{abstract}

\section{Introduction and definition of the ambient space}

Let $\ds \{X_t:M\flecha \overline M\}_{t\in [0,T[}$,  be   a smooth family of immersions of a compact $n$-dimensional manifold (may be with boundary) $M$ into a $(n+1)$-dimensional Riemannian manifold $\overline M$.  By \lq\lq smooth family" we understand that the map $X:M\times [0,T[ \flecha \overline M$ defined by $X(p,t)= X_t(p)$ is smooth. Let $N_t$ be the outward unit normal vector of the immersion $X_t$ and  $H_t$  the trace of the Weingarten map $L_{N_t}$ of $X_t$ associated to $N_t$ (with the convention that $H_t$ is $n$ times the usual mean curvature with the sign which makes positive the mean curvature of a round sphere in $\re^{n+1}$).  By 
$M_t$ we shall denote both the immersion $X_t:M\flecha \overline M$ and the image $X_t(M)$, as well as the Riemannian manifold $(M,g_t)$
 with the metric $g_t$ induced by the immersion. Analogous notation will be used when we have a single immersion $X:M\flecha \overline M$.  We are working with embeddings, then there is no confusion about considering objects on $M$ or on $X_t(M)$.

The {\it volume preserving mean curvature flow} is defined as a solution of the equation
\begin{equation}\label{vpmf}
\parcial{X_t}{t} = (\oH_t- H_t)\ N_t,
\end{equation}
where $\oH_t$ is the averaged mean curvature
$
\oH_t =\ds\frac{ \int_{M_t} H_t dv_{g_t} }{\int_{M_t} dv_{g_t} },
$
being $ dv_{g_t}$  the volume element on $M_t$. This flow decreases the area of $M_t$, but preserves the volume of the domain $\Omega_t$ enclosed by $M_t$ (when such $\Omega_t$  exists).

A nearest relative of \eqref{vpmf} is the mean curvature flow $\parcial{X_t}{t} = - H_t\ N_t$, which has been extensively studied in Euclidean and non Euclidean ambient spaces. The work for the flow \eqref{vpmf} is much shorter and mainly devoted to the evolution of hypersurfaces in the Euclidean space. Special attention has been paid, in the Euclidean or Minkowski ambient spaces,  to the evolution of convex hypersurfaces  (\cite{Hu87} and \cite{An}) and hypersurfaces in the neighborhood of a round sphere (\cite{ES}). The results for the last situation were extended in \cite{AF} to general Riemannian manifolds whose scalar curvature has nondegenerate critical points. Results  both for the evolution of $h$-convex hypersurfaces and hypersurfaces which are near a geodesic sphere were obtained in  \cite{CaMi1} for the Hyperbolic space.  

As far as we know, few other evolutions under \eqref{vpmf} have been studied. This makes the works \cite{Ath1} and  \cite{Ath2} by M. Athanassenas  specially interesting and valuable. Indeed, she deals with two extra difficulties: lack of convexity or sphere closeness assumptions and boundary conditions. In particular, 
she considers the following 

\medskip
{\bf Setting  Euc (\cite{Ath1} and  \cite{Ath2})}.
{\it Let $M$ be a smoothly embedded hypersurface in $\re^{n+1}$ contained  
in the domain
$G = \{x \in \re^{n+1} : 0 \le x^{n+1} \le d;\  d>0\}$,
and with boundary $\partial M \subset \partial G$.  Assume that $M$ is a revolution hypersurface around the axis $z:=x^{n+1}$ generated by the graph of a function $r(z)$, that it intersects $\partial G$ orthogonally at $\partial M$
and encloses (inside $G$) a $(n+1)$-volume $V$.} 
\medskip

Athanassenas' papers study the behavior of the  flow defined by \eqref{vpmf}  with $M$ in the setting Euc as initial condition and the boundary condition that the boundary $\partial M_t$ of the evolving surface $M_t$ remains in $\partial G$ and meets it orthogonally.
 Next we recall her main results:
\begin{itemize}
\item[{\bf A1.}]  The flow  exists for all time
$t > 0$ and converges to the cylinder $C\subset G$ enclosing a volume $V$ under the assumption $n$-volume$(M) \le \fracc{V}{d}$ (cf. \cite{Ath1}).

\item[{\bf A2.}] The singular set (if not empty) is finite and is located along the axis of rotation (cf. \cite{Ath2}).
\end{itemize}

In the present paper, we face the challenge of extending the above statements to hypersurfaces living in a larger universe; in other words, our goal is to export the Euclidean results to more general ambient spaces. But not every Riemannian manifold is well suited to deal with the situation analogous to the setting Euc. Indeed, in order   the concept of revolution hypersurface makes sense, we need that the ambient space also has some kind of rotational symmetry. 
This motivates the introduction of the spaces defined below.

\begin{defi} \lb{AE}
A rotationally symmetric space with respect to an axis $z$ is a manifold $\oM= J \times B^{n}_R$, with $J$ an interval in $\re$ and  $B^n_R$ the Euclidean open ball of radius $R$ in $\re^n$, endowed with a $C^\infty$ metric of the form
\begin{equation}\label{mRSS}
 \overline g:= dr^2 + f^2(r)\  dz^2 + h^2(r) \ g_{S},
\end{equation}
in the cylindrical coordinates $r=|x|$, $z=z$, $u=x/|x|$ of $(z,x)\in \oM$,
where  $g_{S}$ is the standard metric of sectional curvature $1$ on the sphere $S^{n-1}$.
\end{defi}

On such a manifold there is a natural action of $SO(n)$ in the following way. Given a point $q=(z,x)$ and ${\mc R}\in SO(n)$,  we define ${\mc R}q:= (z, {\mc R}x) \equiv \exp_p {\mc R}\exp_p^{-1}q$, where let $p=(z,0)$ and $ \exp_p$ is the exponential map in $(\overline M, \overline g)$. The axis $z$ is invariant under such action, and this is the reason for the name of these spaces. From now on, we shall call the axis $z$ the {\bf axis of rotation}.

Notice that the previous definition includes the three model spaces with constant sectional curvature. In fact,
\begin{itemize}
\item  When $R=\infty$, $J=\re$, $\lambda <0$, $f(r)=\cosh (\sqrt {|\lambda|}\  r)$, $h(r)= \fracc{1}{\sqrt {|\lambda|}} \sinh (\sqrt {|\lambda|}\  r)$, $(\oM,\og)$ is the  Hyperbolic space of sectional curvature $\lambda$. 

\item When $\lambda >0$, $R=\frac{\pi}{2 \sqrt {\lambda}}$, $J=]-\frac{\pi}{2 \sqrt{\lambda}}, \frac{\pi}{2 \sqrt {\lambda}}[$,  $f(r)=\cos (\sqrt {\lambda}\  r)$, $h(r)= \fracc{1}{\sqrt {\lambda}} \sin (\sqrt {\lambda}\  r)$,  $(\oM, \og)$ is the open half sphere of sectional curvature $\lambda$.

\item When $R=\infty$, $J=\re$,  $f(r)=1$, $h(r)= r$,  $(\oM, \og)$ is the Euclidean Space. 
\end{itemize}

According to the previous definition, a natural setting for the extension of the results A1 and A2 mentioned above is:

\medskip
{\bf Setting RSS.} {\it Let $M$ be a smoothly embedded hypersurface of revolution around the axis of rotation in $\oM$, generated by the graph of a function $r(z)$ and contained  
in the domain
$G = \{(r,z,u)\in \oM  : a \le z\le b\}$,  with boundary $\partial M$, which intersects $\partial G$ orthogonally at the boundary and encloses a $(n+1)$-volume V inside $G$.}
\medskip

However, as we pointed out in \cite{CaMi1}, the negative ambient curvature seems to be a friendlier setting for the flow \eqref{vpmf}, so here we still need suitable hypotheses on the negativity of some sectional curvatures. In particular, our main results are valid for the following

\medskip
{\bf Setting RSS2.} {\it Let M be in the setting RSS and assume that either
\begin{itemize}
\item[(a)] the curvatures of $\oM$ satisfy $S_{zi}< 0$ and $S_{ri}\le 0$, or
\item[(b)] $\oM$ is the Euclidean space.
\end{itemize}
Here (cf. notation of section \ref{RSS}), for each $i$, $S_{zi}$ means the sectional curvature of the plane generated by the axis of rotation and a direction $E_i$ orthogonal to the plane $\pi_{G}$ containing the generating curve, and $S_{ri}$ is the sectional curvature given by $E_i$ and the direction in $\pi_{G}$ orthogonal to the axis.}
 \medskip
 
 In this setting,  {\it  flow $M$  by \eqref{vpmf} with the boundary condition that the evolving hypersurface
 \be\lb{bocon}
 \text{ $M_t$ intersects $G$ orthogonally at the boundary for every $t$.}
 \ee}
   We shall prove that

 \begin{itemize}
\item[{\bf B1.}]  The flow 
exists  
while the evolving hypersurface $M_t$ does not touch the axis of rotation (Theorem \ref{t:preLTE}). 
 
\item[{\bf B2.}] As long as the flow exists, the generating curve remains 
a graph over the axis of rotation (Theorem \ref{graphbo}).
 
\item[{\bf B3.}] Throughout the same time, the averaged mean curvature $\oH$ is  positive and double side bounded (Proposition  \ref{t:fboH}).

\item[{\bf B4.}] If the motion $M_t$ of $M$ has a first singularity at time $T$, the singular set is finite and lies on the axis of rotation (Theorem \ref{cdps}).

\item[{\bf B5.}] If  $\vle(M)$ is small enough compared with the volume $V$ enclosed by $M$ inside $G$, then the solution of \eqref{vpmf} is defined for all $t>0$ and converges to a revolution hypersurface of constant mean curvature in $\oM$ (Theorem \ref{fth}).
\end{itemize}

 At this moment, we think that the non positivity of the sectional curvature  is an essential hypothesis for B2, but  we have no concluding counterexample.  The other results may still be true in an ambient space  with positive curvature, but this would require a different kind of argument, because our proofs rely strongly on B2.
 
About the techniques used to prove the above results, this paper follows the ideas introduced in  \cite{AAG} and \cite{Sim} (which study the mean curvature flow for closed rotationally symmetric hypersurfaces) with the clever modifications proposed by M. Athanassenas to deal with the averaged mean curvature, and some methods learned from \cite{ChMa}.

 But these ideas are not enough to work in the setting RSS2, where we actually have to deal with an ambient space non (necessarily) flat and with sectional curvatures non (necessarily) constant. The key to overcome this further difficulty is the mixture of some ideas from the aforementioned works with a bunch of tricky technicalities which spread out all over our paper.
 
 In addition, we had to face an extra and unexpected complication: as far as we are aware, there is a couple of mistakes in \cite{Ath1}, \cite{Ath2} (in the last addend of the first displayed formula of the proof of Proposition 5 on page 63 of \cite{Ath1} and on the last displayed expression on page 10 of \cite{Ath2}); and, at least for the mistake in the computation of \cite{Ath2}, when doing it properly, the arguments that follow fail.
This makes necessary   (and not superfluous by far) the inclusion of the Euclidean space within the setting RSS2. In short, we not only extend the desired results to a wider setting, but we also complete the proof in the Euclidean space. 
 
 The paper is organized as follows. We begin with the study of the main features of the spaces introduced in Definition \ref{AE} (section \ref{RSS}). After some remarks on short time existence and basic computations (section \ref{STE}), we find bounds for the distance to the axis of rotation and for $\oH$ (section \ref{broH}) which will be the key to prove B2 (section \ref{graph}). Section \ref{preLTE} is devoted to the proof of B1, and in section \ref{bfoH} we obtain B3 together with finer bounds for $\oH$. We shall use these finer bounds to prove B4 (section \ref{dis}). Finally, B5 is proved in section \ref{volp}.
 
\bigskip
\noindent {\bf Acknowledgments}\\
{\small 
The first author was supported by a Grant of the {\it Programa Nacional de Formaci\'o  n del Profesorado Universitario} ref: AP2003-3344. Both authors were partially supported  by DGI(Spain) and FEDER Project MTM2007-65852, by the project ACOMP07/088 of {\it Generalitat Valenciana} and by the network: MTM2006-27480-E.}

\section{Rotationally symmetric spaces around an axis}\label{RSS}

To compute the covariant derivatives and curvatures corresponding to the metric \eqref{mRSS}, we shall use the orthonormal frame 
\begin{equation}\label{of}
E_0:= E_r = \parcial{}{r}, \quad E_1:= E_z = \frac{1}{f(r)}\parcial{}{z},\quad  E_i = \frac{1}{h(r)} e_i, \quad i=2, ... , n
\end{equation}
where $\{e_i\}$ is an orthonormal frame of $S^{n-1}$ with its standard  metric, and the dual orthonormal frame
\begin{equation}\label{ocf}
\theta^r = dr, \quad \theta^z = f(r)\ dz, \quad \theta^i = h(r) \delta^i,
\end{equation}
being $\{\delta^i\}$ the dual frame of  $\{e_i\}$ in $S^{n-1}$ with its standard  metric. In these frames, the Cartan connection forms $\omega_a^b$ (defined by $d\theta^b = - \sum_{a=0}^n\omega^b_a \wedge  \theta^a$) are
\begin{equation}\label{cofotg}
\omega^z_r = \frac{f'(r)}{f(r)} \theta^z, \quad \omega^i_r = \frac{h'(r)}{h(r)} \theta^i, \quad \omega^i_z = 0, \quad \omega^i_j = \ ^S\omega^i_j,
\end{equation}
where $^S\omega^i_j$ are the connection $1$-forms  of $S^{n-1}$ with its standard metric  and $'$ denotes the derivative respect to $r$.

From the above expressions, the covariant derivatives of the orthonormal frame (given by ($\oN_XE_a = \sum_{b=0}^{n}\omega_a^b(X) E_b$) are
\begin{align}
\oN_{E_r} E_r= 0, &\qquad \oN_{E_z}E_r = \frac{f'(r)}{f(r)} E_z, \qquad  \oN_{E_i} E_r= \frac{h'(r)}{h(r)} E_i,  \nonumber\\
 \oN_{E_r} E_z= 0, & \qquad \oN_{E_z} E_z= - \frac{f'(r)}{f(r)} E_r, \qquad   \oN_{E_i} E_z=0, \label{code}\\
 \oN_{E_r} E_i= 0, &\qquad \oN_{E_z} E_i= 0, \qquad  \oN_{E_i} E_j= -\frac{h'(r)}{h(r)} \delta_{ij} E_r +  {}^S\omega_j^k (E_i)E_k. \nonumber
\end{align}

Also from \eqref{cofotg} we can compute the curvature forms $\Omega_a^b$ (which satisfy $\Omega_a^b = d\omega_a^b - \sum_{c=0}^n \omega_a^c \wedge \omega_c^b$):
\begin{align} \label{car_om}
\Omega_r^z = \frac{f''}{f}  dr \wedge \theta^z, &\qquad \Omega_r^i= \frac{h''}{h} dr\wedge \theta^i ,  \nonumber\\
 \Omega_z^i = - \frac{h' f'}{hf} \theta^i\wedge \theta^z, & \qquad \Omega_j^i = \ ^S\Omega^i_j - \(\frac{h'}{h}\)^2 \theta^i \wedge \theta^j, 
\end{align}
where $^S\Omega^i_j$ are the curvature $2$-forms  of $S^{n-1}$ with its standard metric.

Let us remark a fact which is necessary to have in mind when doing computations:
$$
\ ^S\omega^i_j (E_k)= \frac1{h} {}^S\omega^i_j (e_k)\qquad \ {}^S\Omega^i_j (E_k,E_\ell)=  \frac1{h^2} {}^S\Omega^i_j (e_k,e_\ell).
$$
This is interesting to recall because $^S\omega^i_j (e_k)$ and $^S\Omega^i_j (e_k,e_\ell)$ give the standard values in $S^{n-1}$ with its standard metric.

\begin{nota}\label{fh0}
In any Riemannian manifold $\mathcal M$, given a point $p\in \mathcal M$ and a geodesic $\Gamma$ through $p$, the eigenvalues of the Weingarten map of the geodesic tube of radius $r$ around $\Gamma$ approach $\infty$ as $1/r$ when $r\to 0$ in the directions tangent to geodesic discs orthogonal to $\Gamma$ (cf. \cite{Gr}). Notice that in the metric  \eqref{mRSS} the Weingarten map of such a tube in the directions tangent to geodesic discs is given by $ \oN_{E_i} E_r$. Then, by formulae \eqref{code}, $\lim_{r\to 0} \fracc{ \frac{h'(r)}{h(r)}}{\frac1{r}}=1 $. 

On the other hand, $h(r)$ is the norm of the vector $e_i$ tangent to the geodesic sphere in the hyperplane $z=$constant with radius $r$, then $h(0):= \lim_{r\to 0} h(r) =0$ with $\lim_{r\to 0} \frac{h(r)}{r}=1$, which combined with the above limit gives $h'(0) = 1$. Moreover, $f(r)$ is the norm of the  vector tangent to the curve $z \mapsto (r_0,z,u_0)$, which is a geodesic when $r_0=0$. Therefore, $f(0)>0$ and constant, and we can take, without loss of generality, $f(0)=1$. 

If we do not have any singularity in our manifold $\oM=J\times B^n_R$ at the points $(z,0)$, the limit of $\nablaa_{E_z} E_z$ must be independent on the direction we approach the point $(z,0)$. But from \eqref{code} we have opposite vectors $\nablaa _{E_z} E_z$  of norms $|\frac{f'(0)}{f(0)}|$ when we approach $(z,0)$ by opposite directions, then we must have $f'(0)=0$.

From all these remarks we have that, if the metric $\overline g$ given by \eqref{mRSS} is a Riemannian metric, then
\begin{equation}\label{signosf'}
f(0)=1, \quad f'(0)=0, \quad h(0)=0, \quad h'(0)=1.
\end{equation}
\end{nota}

Now we use the expression $R_{abcd} = \Omega^a_b(E_c,E_d)$ of the curvature components and \eqref{car_om} to obtain the  sectional curvatures:
\begin{align}
&S_{rz} = R_{rzrz} = - \frac{f''}{f} , 
\qquad S_{ri} = R_{riri}= - \frac{h''}{h}  ,  \nn \\
& S_{zi} = R_{zizi} = - \frac{h' f'}{hf}, 
\qquad S_{ij} = R_{ijij}=  \frac{1 - h'^2}{h^2} . \label{scoM}
\end{align}

In the next two remarks, we translate the curvature conditions of the setting RSS2 in terms of the functions $f$ and $h$ defining $\og$.

\begin{nota} \label{hincreas}
If $S_{ri} \le 0$ the second formula in \eqref{scoM}  implies $h'' \ge0$, then $h'$ is non-decreasing and, since $h'(0)=1$, $h'>0$. Thus $h$ is increasing. 
\end{nota}

\noindent Henceforth, we shall denote by $z(\phi)$ the first positive zero of a function $\phi:\re \flecha \re$. 

\begin{nota}\label{fincreas}
If $S_{zi} \le 0$ {\rm (resp. $S_{zi}< 0$)} 
the third formula in \eqref{scoM}  implies $f'h' (r)\ge 0$ {\rm (resp. $f'h' (r)> 0$) }
for every $r$ and, since $h'(r)>0$ for $r\in [0, z(h')[$,  $f'(r)\ge 0$ {\rm (resp. $f'(r)> 0$) }
for $r$ in this interval. So $f$ is non-decreasing, {\rm (resp. increasing)}
in particular, $f(r) \ge f(0)=1$ {\rm (resp. $f(r) > f(0)=1$)}. 
The property $f'(r)\ge 0$, together with $f'(0)=0$, implies $f''(r)\ge 0$ {\rm (resp. $f''(r)> 0$) }
at some points of $[0, z(h')[$, and $S_{rz} \leq 0$ {\rm (resp. $S_{rz} < 0)$ }
at those points. 

If, moreover,  $S_{ri}\le 0$,  Remark \ref{hincreas} implies $z(h')=\infty$, then the inequalities stated in the last paragraph for an interval now hold  for every $r>0$. 
\end{nota}

\section{Short time existence and some formulae}\label{STE}

Let us begin with a remark on the notation. Along all the paper we use many quantities (functions, vector and tensor fields, ...) depending on the evolved hypersurface $M_t$ of $M$ under the flow \eqref{vpmf}. When they are introduced for the first time,  we write the subindex $_t$ to denote their dependence on $t$, but  later we shall omit it unless it is not clear from the context that we are making reference to the function depending on $t$. Sometimes we also denote the dependence on $t$ by $ (\, .\, , t)$.

Observe that, since $M$ is invariant by the action of $SO(n)$  as a group of  isometries of $\oM$, also $M_t$ will be.  Then the unit normal vector $N_t$ to $M_t$ will be contained in the plane generated by $E_r$ and $E_z$ and can be written as
\be
N= \<N,E_r\> E_r + \<N,E_z\> E_z, \label{N}
\ee
in turn, the unit vector $\t_t$ tangent to the generating curve will be
\be
\frak{t} = -\<N, E_z\> E_r + \<N, E_r\> E_z. \label{t}
\ee

We shall use the coordinates $(r_t,z_t,u_t)$ for $M_t$. Since $M$ is generated by the graph of a function $r(z)$,  for small $t$ a solution of \eqref{vpmf} will still be generated by the graph of a function $r_t(z_t)$. Using this function, the vectors $\frak{t}$ and $N$ can be written as
\be\label{tNf}
\frak{t}= \frac{1}{\sqrt{\dot{r}^2 + f^2}}(\dot{r} E_r + f E_z), \quad N= \frac{1}{\sqrt{\dot{r}^2 + f^2}} (f E_r - \dot{r} E_z),
\ee
where $\dot{r}$ denotes the derivative of $r_t$ with respect to $z_t$ . 

A convenient orthonormal frame of $M_t$ is given by $\frak{t}, E_2, ... , E_n$. If we compute the mean curvature of $M_t$ using this frame, and denote $\tilde{\frak{t}}= \sqrt{\dot{r}^2 + f^2} \  \frak{t}$, $\tilde{N}= \sqrt{\dot{r}^2 + f^2} \  N$, we get
\be\label{Hf1} H = k_1 + (n-1) k_2, \ee
where
\begin{align} 
k_1 = - \frac{1}{(\sqrt{\dot{r}^2 + f^2})^3} \<\oN_{\tilde{\t}} \tilde{\t}, \tilde{N}\>
&= \frac{1}{\sqrt{\dot{r}^2 + f^2}} \left( \frac{-\ddot{r} f + \dot{r}^2 f'}{\dot{r}^2 + f^2} + f'\right) \label{k1f1}\\
&= \frac{1}{\sqrt{\dot{r}^2 + f^2}} \left(- \frac{d}{dz}\arctan\left(\frac{\dot{r}}{f}\right) + f'\right)\label{k1f2}
\end{align}
is the normal curvature of $M_t$ in the direction of $\frak{t}$, and
\be\label{k2f} k_2 =\frac{1}{\sqrt{\dot{r}^2 + f^2}}  \<\oN_{E_2}\tilde{N} , E_2\> = \frac{f}{\sqrt{\dot{r}^2 + f^2}} \frac{h'}{h}\ee
is the nomal curvature of $M_t$ in the direction of $E_i,\ i=2, ..., n$.

It is well known (cf. \cite{Eck}) that, up to tangential diffeomorphisms, equation \eqref{vpmf} is equivalent to 
\be\label{vpmft}
\<\parcial {X_t}{t}, N_t\> = \oH - H.
\ee
In this flow the variable $z$ does not change with time, and formulae \eqref{tNf}, \eqref{Hf1}, \eqref{k1f1} and \eqref{k2f} are true just by this change. Using them,  equation \eqref{vpmft} becomes
\be\label{tafl}
\parcial{r}{t} = \frac{\ddot{r}}{\dot{r}^2 + f^2} - \frac{f'(r)}{f(r)} \(1 + \frac{\dot{r}^2}{\dot{r}^2 + f^2}\) - (n-1)\frac{h'(r)}{h(r)} + \oH \sqrt{1 + \frac{\dot{r}^2}{f^2}}.
\ee
Replacing $\oH$ in \eqref{tafl} by any $C^{1, \alpha/2}$ real valued function $\psi$ such that $\psi(0) = \oH(0)$, we obtain a parabolic equation which, at least for small $t$, has a unique solution satisfying $\dot{r}(a) =\dot{r}(b)=0$. Now, using a routine fixed point argument (cf. \cite{Mc3}), we can establish short time existence also for \eqref{tafl} with the same boundary conditions.

\section{Upper bounds for $r$ and  rough bounds for $\oH$}\label{broH}

In this section, we shall prove  that if $M_t$ is a maximal solution  of \eqref{vpmf} with initial condition in the setting RSS and satisfying the boundary condition \eqref{bocon}, then $r_t$ has a finite upper bound and $\oH$ is positive and bounded from above. Under additional curvature assumptions, which include both situations of the setting RSS2, we shall show that $\oH$ is actually bounded away from zero.

The aforementioned bounds will be the key to prove the conservation of the property of being a graph for the generating curve (cf. section \ref{graph}). This and the results of section \ref{preLTE} will help us to obtain more accurate bounds for $\oH$ (cf. section \ref{bfoH}).

Let us define the functions $\beta$ and $\delta$ and the numbers $r_2 > r_1>0$ by 
\begin{align}
 \beta(r)=\int_0^r f(r) h(r)^{n-1} dr, \qquad & r_1=\beta^{-1}\left(\fracc{V}{(b-a) \sigma}\right),\label{defr1}\\
  \delta(r)=\int_0^r  h(r)^{n-1} dr, \qquad & r_2=\delta^{-1}\left(\fracc{\vle(M)}{ \sigma} + \delta(r_1)\right),\label{defr2}
\end{align}
where $\sigma$ denotes the volume of $S^{n-1}$ with its standard metric.
Let us remark that the functions $\beta $ and $\delta$ are increasing, thus they have inverse and $r_2>r_1$ as claimed. Moreover, $r_1$ and $r_2$ depend only on $b-a$, $V$,  $\vle(M)$ and the metric $\overline g$ of the ambient space $\oM$.

 \begin{prop}\label{bor}
 If $M$ is in the setting RSS and $[0,T[$ is the maximal time interval  where the flow  \eqref{vpmf} satisfying \eqref{bocon} is defined, then  $r_t < r_2$ for every $t\in[0,T[$.
\end{prop}
\begin{demo} 
From the expression  \eqref{mRSS} of the metric of $\oM$ it follows that its volume element is 
\be\label{muaf}
\overline{\mu}  = f(r) h(r)^{n-1} dr \wedge dz \wedge \ ^S\mu ,
\ee
 being  $^S\mu$ the volume element on $S^{n-1}$ defined by its standard metric. 

First, let us suppose that the generating curve is a graph of a function for all time $t \in [0,T[$. Using expression \eqref{muaf},  the volume $V$ enclosed by $M_t$ can be computed as
\be
V= \sigma \int_a^b \int_0^{r(z)} f(r) h(r)^{n-1} dr \ dz,  \lb{fVrz}
\ee
and  the definition of $r_1$ gives
\be
V= \sigma \int_a^b \int_0^{r_1} f(r) h(r)^{n-1} dr \ dz.  \label{fVr1}
\ee
Comparing these formulae, we deduce that $\min_{z\in[a,b]}r(z) =: r_m \le r_1 \le r_M:= \max_{z\in[a,b]}r(z)$ and one of the inequalities is strict if and only if the other is.

The induced volume form $\mu_t$ on $M_t$ is defined by the contraction
\be\label{muf} \mu_t = \imath_{N_t} \overline{\mu} =  \sqrt{\dot{r}^2(z) + f^2(r(z))} \  h(r(z))^{n-1} dz \wedge \ ^S\mu,
\ee
where we have used \eqref{ocf}, \eqref{tNf} and \eqref{muaf}. Hence
\begin{align}
\vle(M_{t}) &=   \sigma \int_a^b  \sqrt{\dot{r}^2(z) + f^2(r(z))} \  h(r(z))^{n-1} dz >  \sigma \int_a^b  |\dot{r}|\  h(r(z))^{n-1} dz \nonumber\\
&\ge  \sigma \int_{r_m}^{r_M}     h(r)^{n-1} dr =   \sigma \left(\int_{r_m}^{r_1}     h(r)^{n-1} dr + \int_{r_1}^{r_M}     h(r)^{n-1} dr\right) \label{vle>+}.
\end{align}
From this inequality we obtain 
\be
\vle(M_t) > \sigma  \int_{r_1}^{r_M}     h(r)^{n-1} dr = \sigma (\delta(r_M) - \delta(r_1)).
\ee
Then, having into account that $\vle(M_t)< \vle(M)$, we reach
$$\delta(r_M) < \fracc{\vle(M)}{ \sigma} + \delta(r_1),$$
from which the proposition follows.

For the general case, we can parametrize the generating curve as $(r(s),z(s))$, $s\in[0,c]$, then formula \eqref{muf} becomes $ \mu_t =  \sqrt{\dot{r}^2(s) + f^2(r(s)) \dot{z}^2(s)} \  h(r(s))^{n-1} ds \wedge \ ^S\mu,$ and we have to distinguish betwenn the cases $r_M \le r_1$ (then all is done) and $r_M\ge r_1$. In the last case, \eqref{vle>+} still works just changing $dz$ by $ds$, $[a,b]$ by $[0,c]$, the functions of $z$ by the functions of $s$ and $f^2$ by $f^2 \dot{z}^2$. After that, the proof finishes as in the above case.
\end{demo}

Next, the goal is to bound the averaged mean curvature $\oH_t$.

\begin{prop}\label{p:uboH}  Let $M_t$ be the solution of  \eqref{vpmf}  with initial condition $M$ in the setting RSS and satisfying \eqref{bocon}.
For every $t$ such that $$0<\rho \le r_t \le z_0:= \min\{z(h'), z((n-1) h' f + f' h)\} \quad \text{ for some fixed } \rho,$$ 
and the generating curve of $M_t$ is a graph, there is a constant $h_2(V, \overline g,n,b-a,\vle(M),\rho)> 0$   such that $0\le \oH \le h_2$.
\end{prop}
\noindent \begin{demo} Before starting the proof, let us remark that
\begin{itemize}
\item By \eqref{signosf'},  we have $ \min\{z(h'), z((n-1) h' f + f' h)\}>0$.
\item The hypothesis on the variation of $r_t$ together with Proposition \ref{bor} imply that  $r_t \le \frak{r} = \min\{r_2,z_0\}$.
\end{itemize}

From \eqref{Hf1}, \eqref{k1f2} and \eqref{k2f} we can write
\be
\oH_t = \frac1{\vle(M_t)} \int_{M_t}  \frac{1}{\sqrt{\dot{r}^2 + f^2}} \left(- \frac{d}{dz}\arctan\left(\frac{\dot{r}}{f}\right) + (n-1) \frac{h'}{h} f +  f' \right) \mu_t = I_1 + I_2,
\ee
where we have (by \eqref{muf}) the following expressions for $I_1$ and $I_2$:
\begin{align}
I_1 &= \frac{\sigma}{\vle(M_t)}  \int_a^b   - \frac{d}{dz}\arctan\left(\frac{\dot{r}}{f} \right) h(r(z))^{n-1} \ dz 
\quad \text{ and }  \label{I1} \\
I_2&= \frac{1}{\vle(M_t)}  \int_{M_t}  \frac{1}{\sqrt{\dot{r}^2 + f^2}}  \frac{(n-1) h' f + h f'}{h}  \mu_t.  \label{I2}
\end{align}

Integrating by parts and realizing that the condition in the boundary gives $\dot{r}(b)=\dot{r}(a)=0$, we get
\be
I_1=\frac{\sigma}{\vle(M_t)} \int_a^b \arctan\left(\frac{\dot{r}}{f} \right)(h^{n-1})' \  \dot{r} \  dz \ge 0 \label{I1pos}
\ee
because $\arctan\left(\fracc{\dot{r}}{f}\right)  \dot{r} \ge 0$. Since $\arctan\left(\fracc{\dot{r}}{f} \right) \dot{r} \le \fracc{\pi}{2} \sqrt{\dot{r}^2 + f^2}$, 
\begin{align}
I_1  &  \le \frac{\sigma}{\vle(M_t)} \frac\pi2 \int_a^b \sqrt{\dot{r}^2 + f^2}\ (h^{n-1})' \ dz \lb{I1pos2} \\
&= \frac{\pi }{2 \vle(M_t)} \int_{S^{n-1}} \int_a^b \frac{(h^{n-1})'}{h^{n-1}} h^{n-1} \ \sqrt{\dot{r}^2 + f^2} \ dz \ ^S\mu \nonumber\\
& = \frac{\pi}{2 \vle(M_t)} \int_{M_t} \frac{(h^{n-1})'}{h^{n-1}} \mu_t   \le  \frac{\pi }{2 } \max_{r\in[\rho,\frak{r}]} \frac{(h^{n-1})'}{h^{n-1}}(r)<\infty. \label{I1ub}
\end{align}

From \eqref{muf} and \eqref{I2}  we have
\begin{align}\label{I2pos}
I_2 &= \frac{\sigma}{\vle(M_t)} \int_a^b \((n-1) h' f + f' h\) h^{n-2} dz \ge 0.
\end{align}
On the other hand, \eqref{I2}  also implies that
\be\label{I2ub}
I_2 \le  \frac{1}{\vle(M_t)} \int_{M_t} \frac1f \frac{(n-1) h' f + f' h}{h} \mu_t \le  \max_{r\in[\rho,\frak{r}]}   \frac{(n-1) h' f + f' h}{f h}(r) <\infty.
\ee
In conclusion, the existence of the finite upper bound $h_2$ follows from \eqref{I1ub} and \eqref{I2ub}. Moreover, the non-negativity of $\oH$ is due to \eqref{I1pos} and \eqref{I2pos}, which in turn are true thanks to the assumption $r_t \leq z_0$.
\end{demo}

\begin{coro} \label{c:uboH} Let $M_t$ be the solution of  \eqref{vpmf}  with initial condition $M$ in the setting RSS and satisfying \eqref{bocon}. Let us suppose that  the sectional curvatures $S_{zi}$ and $S_{ri}$ of $\oM$ are non-positive. 
For every $t$ such that $0<\rho \le r_t$  for some fixed $\rho$, and the generating curve of $M_t$ is a graph, there are constants $h_i(V, \overline g,n,b-a,\vle(M),\rho)> 0$, $i=1,2$   such that $h_1\le \oH \le h_2$.
\end{coro}
\begin{demo}
 From Remark \ref{fincreas} it follows that the hypotheses on the curvature of $\oM$  imply $z_0=\infty$. Then the upper bound is an immediate consequence of Proposition  \ref{p:uboH}. Remarks \ref{hincreas} and \ref{fincreas} tell us that the hypotheses on the curvature imply that $h$, $h'$ are non decreasing, $f\ge 1$ and $f'h\ge 0$, therefore, from \eqref{I2pos} we obtain
\begin{align}
I_2 &\ge   \frac{\sigma (b-a)}{\vle(M_t)} \min_{r\in[\rho,r_2]}   \((n-1) h' f + f' h\) h^{n-2}(r) \ge  \frac{\sigma (b-a)}{\vle(M)} (n-1) \(h' h^{n-2}\) (\rho) > 0,
\end{align}
which gives the positive lower bound for $\oH$.
\end{demo}

\section{Preserving the property of being a graph for the generating curve}\label{graph}

This section is devoted to prove that, for $M$ in the setting RSS2, the evolving hypersurface remains  a revolution hypersurface generated by a smooth graph. As we pointed out above, $M_t$ is always a revolution hypersurface. Then the aim is to show that the generating curve remains  a graph over the axis of rotation for all the time.

In this section we continue with the notation used in sections 2 and 3. Moreover, by $\a$ we represent the second fundamental form of $M_t$, with the sign convention $\alpha(X,Y) = \<\nablaa_XN,Y\>$, and $L$ denotes its Weingarten map (defined by $LX = \nablaa_X N$). $\nabla$ means the intrinsic covariant derivative on $M_t$.

Recall that the generating curve is a graph if and only if $u:=\< N, E_r\> >0$, which is equivalent to say $1< v = \fracc{1}{u}< \infty$. Therefore, our goal is to obtain an upper bound for $v$. To achieve this,  we need  the evolution equation for $v$.

\begin{lema} Under \eqref{vpmf}, $v= \<N,E_r\>^{-1}$ evolves as
\begin{align}
\parcial{v}{t} 
=& \Delta v - \frac2v |dv|^2- \( \frac1{v} \frac{f'}{f} - k_1 \)^2  v + \oH \frac{f'}{f}(1- v^2)   \nonumber\\
 & - \(\frac{f''}{f} - 2 \(\frac{f'}{f}\)^2 + (n-1)\( \(\frac{h'}{h}\)' - \frac{h'f'}{hf}\)\)   \(1- \frac1{v^2} \)  v.  \label{dvtn}
\end {align}
\end{lema}
\begin{demo} First we compute $\Delta u$. To do so, we shall use an orthonormal frame of $M$ of the form $\t, E_2, ... , E_n$. It follows from \eqref{code}, \eqref{N}  and \eqref{t} (compare with  \eqref{k2f})  that 
\begin{align}
\oN_{E_i} N &= k_2 E_i = u\  \frac{h'}{h} E_i,   \label{5.2}\\
E_i  u &=0,  \qquad \qquad \nabla_\t \t =0, \label{ei<ntt}\\
  \oN_\t E_r &= u\ \frac{f'}{f} E_z,  \qquad \oN_\t N = k_1 \t . \label{5.4}\end{align}
From \eqref{ei<ntt} we have 
\be\label{Deltau}
\Delta u = \t \t u - \sum_{i=2}^n (\nabla_{E_i} E_i) u,
\ee

\noindent and, using  formulae \eqref{5.4}, \eqref{code}, 
\eqref{t} and $\<N,E_z\>^2 = 1 -u^2$, we obtain
\begin{equation*}
\t u = \(u \frac{f'}{f} -k_1\) \sqrt{1-u^2},
\end{equation*}
\begin{align}
 \t \t u&=   \(u \frac{f'}{f} -k_1\) \frac{f'}{f} (1-u^2)  - \(u \frac{f'}{f} -k_1\)^2 u \nonumber \\
  &  \quad- u \( \frac{f'}{f}\)' (1-u^2) - \t k_1\ \sqrt{1-u^2}, \label{ttu}
\end{align}
\begin{align}
(\nabla_{E_i} E_i) u&= \<\nabla_{E_i} E_i, \t\> \t u = \frac{h'}{h}  (1-u^2)  \(u \frac{f'}{f} -k_1\).\label{neeu}
\end{align}
Joining \eqref{Deltau}, \eqref{ttu} and \eqref{neeu}, we reach the desired expression for $\Delta u$:
\begin{align} \label{5.8}
\Delta u = &  \(u \frac{f'}{f} -k_1\)  \frac{f'}{f} (1-u^2) -\sqrt{1-u^2}\  \t k_1 -  \(u \frac{f'}{f} -k_1\)^2 u \nonumber \\
&  - u \(\frac{f'}{f}\)' (1-u^2) - (n-1)  \(u \frac{f'}{f} -k_1\)  \frac {h'}{h} (1-u^2) .
\end{align}

From the standard evolution formula $\fracc{\oN N}{\partial t} = \grad H$ for \eqref{vpmf}, using again \eqref{code} and \eqref{t}, we get
\begin{align}
\parcial{}{t} u &= \<\grad H, E_r\> + (\oH-H)  \<N,E_z\>^2 \frac{f'}{f} \nonumber
\end{align}
but
\begin{align} \label{5.9}
 \<\grad H, E_r\>  &= (\t k_1+ (n-1) \t k_2 )  \<\t, E_r\> \nonumber \\
&= - \sqrt{1-u^2} \ \t k_1 + (n-1) u (1-u^2) \(\frac{ h'}{h} \)'   \nonumber \\
& \quad - (n-1) (1-u^2) \(u \frac{f'}{f} - k_1\) \frac{h'}{h}. 
\end{align}
Then, combining \eqref{5.8} with \eqref{5.9} and rearranging terms, we arrive at the formula
\begin{align}
\parcial{u}{t} = &  \Delta u +  \(u \frac{f'}{f} - k_1 \)^2 u + u (1-u^2) \(\frac{f''}{f} - 2 \(\frac{f'}{f}\)^2 + (n-1)\( \(\frac{h'}{h}\)' - \frac{h'f'}{hf}\)\)\nonumber \\
&  +  \oH \frac{f'}{f} (1-u^2)  \label{dutn}
\end{align}

Now the definition of $v$ leads to
\be\label{vfu}
dv = - \frac1{u^2} du, \quad \parcial{v}{t} = - \frac1{u^2} \parcial{u}{t}, \quad \Delta v = -\frac1{u^2} \Delta u + \frac{2}{u^3} |du|^2.
\ee
Finally, \eqref{dvtn} follows from \eqref{dutn} and \eqref{vfu}.
\end{demo}

\begin{teor}\label{graphbo} 
Let $M_t$ be the solution of  \eqref{vpmf}  defined on a maximal time interval $[0,T[$, with initial condition $M$ in the setting RSS2 and satisfying \eqref{bocon}. The generating curve of the solution $M_t$ of \eqref{vpmf} remains a graph over the axis of revolution for every $t\in[0,T[$. 
\end{teor}
\begin{demo} The case of $\oM$ being the Euclidean space was proved in \cite{Ath1}. Hence we shall suppose that $S_{zi}<0$ and $S_{ri}\le 0$.

For every $t_0\in[0,T[$, let $\rho(t_0) = \min_{(x,t)\in M\times [0,t_0]} r_t(x) >0$. Let $t_1\in[0,t_0]$ be any time such that 
the generating curve of $M_t$ is a graph for every $t\in[0,t_1]$. Then we can apply Corollary \ref{c:uboH} to conclude that 
\be\lb{h1Hh2}
0 < h_1 \le \oH_t \le h_2 \qquad \text{ for } t\in [0,t_1],
\ee
where $h_i := h_i(V, \overline g,n,b-a,\vle(M),\rho(t_0))$.

While the generating curve is a graph, one has $v\ge 1$ and so 
\be\lb{oH<oHv}
\(1-\fracc1{v^2}\) v \le v \quad \text { and } \quad  \oH \frac{f'}{f}  \le  \oH \frac{f'}{f}  v.
\ee 

Putting the above inequalities in  \eqref{dvtn} and forgetting about the negative addends, we reach the inequality:
\begin{align} 
\parcial{v}{t} &\le \Delta v - \frac2v |dv|^2 - \( \frac1v \frac{f'}{f} -k_1\)^2 v - K v^2 + C v\le \Delta v  - K v^2  + C \, v, \label{indvt*}
\end{align}
where
\be\lb{K}
K  = K(V, \og,n,b-a,\vle(M),\rho(t_0)) = h_1 \min_{r\in[ \rho(t_0), r_2]}\frac{f'}{f} ,
\ee
\begin{align} 
C & =C(V, \og,n,b-a,\vle(M),\rho(t_0)) \nn \\
&= \max_{r\in[ \rho(t_0), r_2]} \left|
   \frac{f''}{f} - 2\(\frac{f'}{f}\)^2 + (n-1) \( \left(\frac{h'}{h}\right)' - \frac{f'h'}{fh} \) \right|   + h_2  \max_{r\in[ \rho(t_0), r_2]}\frac{f'}{f}.   \label{C}
\end{align}

Let $(x_2,t_2)\in M_{t_2} \times [0,t_1]$ be the space-time point where $v$ attains its maximum for $t\le t_1$. If $v_2:=v(x_2,t_2) >1$ and $t_2\ne 0$, then $x_2$ is an interior point of $M_{t_2}$ (because at the boundary $v=1$) and it has to satisfy $\parcial{v}{t} \ge 0$  as well as $\Delta v \le 0$. Using these inequalities in \eqref{indvt*}, we can write
$0 \le - K v_2^2  + C \, v_2 = v_2 (-K v_2 + C)$, and so 
$$v_2 \le C_1= \max\{\max_{x\in M} v(x), C/K\}.$$
Notice that $C_1$ does not depend on $t_1$, thus $v$ cannot become $\infty$ in the interval $[0,t_0]$.  Since $t_0$ is arbitrary, the same conclusion is true for $[0,T[$. 
\end{demo}

\section{The first singularities of the motion are produced at the axis of revolution}\lb{preLTE}

This section is devoted to the proof of the following result, which assures long time existence unless the evolving hypersurface reaches the axis of rotation.

\begin{teor}\lb{t:preLTE}
Let $M_t$ be the maximal solution of  \eqref{vpmf}, defined on  $[0,T[$, with initial condition $M$ in the setting RSS2 and satisfying \eqref{bocon}. Then either $$T=\infty \qquad \text{ or } \qquad \lim_{t\to T} \min_{x\in M_t} r_t(x) = 0.$$ 
\end{teor}
\begin{demo}
First, let us recall that $r_t>0$ for every $t\in [0,T[$. Hence the continuous function $t\mapsto \beta(t):= \min_{x\in M_t} r_t(x) $ satisfies $\beta(t)>0$ for every $t\in [0,T[$. If $\lim_{t\to T}\beta(t)>0$, then there is a $\rho>0$ such that $\beta(t) \ge \rho$ for every $t\in [0,T[$. We shall prove that, if this is the case, then  the solution of the flow can be prolonged after $T$, which is a contradiction.

The key  for proving that the solution can be extended if $r_t\ge \rho > 0$ is to show that, under this condition,  $|L|^2$ is bounded. To do so, we need the evolution equation of $|L|^2$ under \eqref{vpmf} for hypersurfaces in arbitrary Riemannian manifolds. Previously, it is necessary to know the evolution equations for the metric $g_t$ of $M_t$, its dual metric $g_t^\flat$ and the second fundamental form $\alpha_t$:
 \be
   \parcial{g}{t} = 2 (\overline H-  H)  \a, \qquad  \parcial{g^\flat}{t} = - 2 (\overline H-  H)  \a^\flat, \nonumber 
\ee
 \begin{align}
\parcial{\a}{t} &= \Delta \a + \(\oH - 2 \ H \) \a L - \oH \,\oR_{\ \cdot\ N\ \cdot\ N} + (|L|^2 + \oRic(N,N))  \a  \nn
\\ & \quad -  \oRic(\cdot\  ,\ L\  \cdot\ ) \nonumber
-  \oRic(L \ \cdot\  ,\cdot\ ) +  \oR(N,  L \  \cdot\  , N,\ \cdot ) \nonumber 
  + \oR(N,  \  \cdot\  , N, L \  \cdot ) 
\\ &  \quad + 2 RiL - \overline\nabla_\cdot \oRic(\cdot, N) - \overline\delta\oR_{\ \cdot\  \cdot\  N},  \nonumber 
\end{align}
where 
$RiL(Z,Y)= \sum_{i=1}^n  \oR_{e_i,Z,Le_i, Y}$  and $ \overline\delta\oR= 
\sum_{i=1}^n  \overline \nabla_{e_i}(\oR)_{e_i}$ for some local orthonormal frame $\{e_i\}$ of $M_t$ and
  $\a L$ is defined by $\a L(X,Y) = \a(L X, Y)$. 

From these equations  we obtain the evolution formula for the squared-norm of the shape operator $|L|^2$:
\begin{align}\lb{Acu}
 \parcial{|L|^2}{t} &= \Delta |L|^2 - 2|\nabla L|^2 + 2|L|^2\(|L|^2 +\oRic(N,N)\)   \nn
 \\ & \quad -2 \oH \( \tr L^3 + \sum_j \oR_{Le_j N e_j N} \) -4 \sum_{i,j} \(\oR_{Le_i e_j Le_i e_j}-  \oR_{Le_i e_j e_i Le_j}\)  \nn
 \\ & \quad -2 \sum_{i,j} \(\nablaa_{Le_i} \oR_{N e_j e_i e_j} + \nablaa_{e_j} \oR_{Ne_i Le_i e_j}\).
 \end{align}
  
 In our case, using the orthonormal local frame $N,\t=E_1, E_2, ..., E_n$, we get
\begin{align}\lb{evL2}
 \parcial{|L|^2}{t} &= \Delta |L|^2 - 2|\nabla L|^2 + 2|L|^4- 2 \oH  \tr L^3 +  2 |L|^2 \(\mc{T} + (n-1) \mc{J}\)  \nn \\
 & \quad  - 2 \oH \(k_1 \mc{T} + (n-1) k_2 \mc{J} \) - 4 (k_1- k_2)^2 (n-1) \mc{Y} - 2 \< \a  , \tilde{\delta} \oR_N \>,
 \end{align}
where $\tilde{\delta} \oR_N(X,Y):= \sum_i\(\nablaa_X \oR_{N E_i Y E_i} + \nablaa_ {E_i}\oR_{N Y X E_i} \)$, $ \mc{T}$, $\mc{J} $ and $\mc{Y}$ are the sectional curvatures of the planes generated by $\{\t, N\}$, $\{E_i, N\}$ and $\{E_i, \t\}$, respectively.  

We cannot deduce any bound for $|L|$ directly  from \eqref{evL2}. Then, as in \cite{Ath1}, we shall study the evolution of $g = (\varphi\circ v) |L|^2$, where $\varphi$ is a real function to be determined later. The evolution of $g$ is given by those of  $\varphi$ and $|L|^2$ according to the formula:
\begin{align}
&\(\parcial{}{t} -  \Delta\) g = |L|^2 \(\parcial{}{t} - \Delta\) \varphi +  \varphi \(\parcial{}{t} - \Delta\) |L|^2 - 2 \< d \varphi, d |L|^2\> \nn \\
& \qquad = |L|^2 \varphi' \(\parcial{}{t} - \Delta\) v - |L|^2 \varphi''\ |dv|^2 +  \varphi \(\parcial{}{t} - \Delta\) |L|^2 - 2 \< d \varphi, d |L|^2\>. \lb{feg}
\end{align}

For the last addend in \eqref{feg}, we use the following inequality (cf. \cite{EcHu})
\be\lb{EcHu}
- 2 \< d\varphi , d|L|^2 \> \le - \frac1{\varphi} \<d g, d\varphi\> + 2 \varphi |d|L|^2| + \frac3{2 \varphi} |L|^2 |d\varphi|^2
\ee
and Kato's inequality
\be\label{Kat}
| d|L| |^2 \le |\nabla L|^2.
\ee
 If $\varphi'>0$, from \eqref{indvt*} to \eqref{C} and \eqref{evL2} to \eqref{Kat}, we conclude that there are constants $K  = K(V, \og,n,b-a,\vle(M),\rho)$ and
$C  =C(V, \og,n,b-a,\vle(M),\rho)$ such that
\begin{align}
\(\parcial{}{t} -  \Delta\) g
\le & \ |L|^2 \varphi' \( - \frac{2}{v}  |d v|^2   - \( \frac{f'}{v f} - k_1\)^2 v - K v^2 + C v\)    - |L|^2 \varphi'' |dv|^2\nonumber\\
&+ 2 \varphi |L|^4 - 2 \varphi \oH \tr L^3 + 2 \varphi |L|^2 \(\mc{T} + (n-1) \mc{J}\)  \nn \\
 &   -2 \varphi \oH \(k_1 \mc{T} + (n-1) k_2 \mc{J} \) - 4 \varphi (k_1- k_2)^2 (n-1) \mc{Y} \nn
 \\ &  -2 \< \a  , \tilde{\delta}\oR_N \>  \varphi 
- \frac1{\varphi} \< dg, d\varphi \>  + \frac3{2 \varphi} |L|^2 |d\varphi|^2.  \lb{evolg1}
\end{align}

Next, let us bound and/or rearrange the different addends in \eqref{evolg1}.
\be\lb{bdphi2}
 - \frac{2}{v} |L|^2 \varphi'  |dv|^2  - |L|^2 \varphi'' |dv|^2 + \frac3{2 \varphi} |L|^2 |d\varphi|^2 = \(- \frac{2}{v \varphi'}  - \frac{\varphi''}{{(\varphi')}^2} +\frac3{2 \varphi}\) |L|^2 |d\varphi|^2.
\ee
\begin{align}
 2 \varphi |L|^4 - |L|^2 \varphi' v &\( \frac{f'}{v f} - k_1\)^2 = g \frac{\varphi'}{\varphi} v \(2 \frac{k_1 f'}{v f} - \frac{(f')^2}{v^2 f^2} -|L|^2 + (n-1)k_2^2\) + 2 \varphi |L|^4 \nn \\
& = \(- \frac{\varphi'  v}{\varphi^2} + \frac2\varphi \)g^2 + (n-1) k_2^2 \frac{\varphi'  v}{\varphi} g + 2 k_1 \frac{f' \varphi'}{f \varphi} g -  \(\frac{f'}{f}\)^2 \frac{\varphi'}{v \varphi} g.\lb{b-g2}
\end{align}
Using Young's inequality $x y \le \varepsilon x^2 + \frac1{4 \varepsilon} y^2$ with $x=g$, $y = \frac{2 f'  \varphi'  |k_1|}{f \varphi}$ and $\varepsilon = k/2$, for some $k>0$ that we shall determine later, we have
\begin{align}
 2 k_1 \frac{f' \varphi'}{f \varphi} g &\le   \frac{k}{2} g^2 + \frac{1}{2 k} \(\frac{2  k_1 f'  \varphi' }{f \varphi}\)^2 \le \frac{k}{2} g^2 + \frac{2}{k \varphi} \(\frac{\varphi' f'}{\varphi f}\)^2 g. \label{bf'fLH}
\end{align}
If $\varphi$ is bounded from below and $\varphi'$ is bounded from above by some positive constants, there is $0<K_1 = C \sup \( \fracc{\varphi'}{\varphi} v\)$ satisfying
\begin{align}
|L|^2 \varphi' ( - K v^2 + C v)\le &K_1 g. \lb{K1g}
\end{align}
Using $|\tr L^3| \le |L|^3$ and Young's inequality  with $x=|L|$, $y = 2  \oH$ and $\varepsilon = k \varphi$:
\begin{align}
- 2 \varphi  \oH  \tr L^3 \le 2 \varphi \oH |L|^3 = 2  \oH |L| g \le \(k \varphi |L|^2 + \frac1{4 k \varphi} 4 \oH^2 \) g = k g^2 + \frac{\oH^2}{k \varphi} g. \lb{boHtrL}
\end{align}
To bound the rest of the addends in \eqref{evolg1}, we recall that, since $\rho\le r \le r_2$, the curvatures $\mc{T}$, $\mc{J}$, $\mc{Y}$ and $|\tilde{\delta}\oR_N|$ are bounded; moreover, as
$$k_1 \le |L|, \qquad k_2 \le |L| \qquad \text{and} \qquad  (k_1 - k_2)^2 \le 3 |L|^2,$$
there are positive constants $C_1, C_2, C_3$ and $C_4$ such that
\begin{align}2 \varphi |L|^2 &\(\mc{T} + (n-1) \mc{J}\)   +2 \varphi \oH \(k_1 \mc{T} + (n-1) k_2 \mc{J} \) \nn \\
 &  - 4 \varphi (k_1- k_2)^2 (n-1) \mc{Y} 
   -2 \< \a  , \tilde{\delta}\oR_N \>  \varphi \nn \\
   & \le C_1 g + C_2 \sqrt{\varphi} \sqrt{g} + C_3 g + C_4 \sqrt{\varphi} \sqrt{g}. \lb{brest}
 \end{align}

Now, in order that $\varphi'\circ v$ and $\varphi\circ v$ be double  side bounded, and $\varphi > 1$ (conditions that we have used before), we choose the function $\varphi$ and the constant $k$ to be
\be\label{fi}
\varphi(v) := \frac{v^2}{1-k v^2}, \qquad k:= \frac1{2\max{v^2}}.
\ee
With this election, $\varphi$ also satisfies 
\be\lb{otrask}
\frac3{2\varphi} - \frac{\varphi''}{\varphi'^2}-\frac2{\varphi'v}  =  \frac{v^2 - \varphi}{2 \varphi^2} <0 \quad \text{ and } \quad \frac2{\varphi} - \frac{\varphi'}{\varphi^2} v + k+\frac{k}{2} = -\frac{k}{2} .
\ee
Thanks to the expressions from  \eqref{evolg1}  to   \eqref{otrask}, we reach 
\be\lb{evolg3}
\(\parcial{}{t} -  \Delta\) g \le -\frac{k}{2} g^2 +
K_2 g + K_3 \sqrt{g}  - \frac1{\varphi} \< dg , d\varphi\>,
\ee
for two positive constants $K_2$ and $K_3$.

Let $t_1\in[0,T[$. Let us suppose that $(x_0,t_0)$  is the point where $g$ attains its maximum for $t\le t_1$. We further assume that $x_0$ is an interior point of $M_{t_0}$  and $t_0\ne 0$.  Observe that this assumption can be made without loss of generality: 
\begin{itemize}
\item[i.] If $z(x_0)=b$, then we can use the map $\imath_1:(r,z,u) \mapsto (r, 2 b - z, u)$, which is an  isometry of  $\oM$, to extend the solution symmetrically onto $[a, 2 b - a]$. In fact, by applying $\imath_1$ to $M_t$, we construct a new solution $M_t \cup \imath_1(M_t)$ of  \eqref{vpmf}, which is orthogonal to $z=a$ and $z=2 b-a$ at the boundary,  such that $x_0$ still is a maximum of $g$, but it is not at the spatial boundary.
\item[ii.] If $z(x_0)=a$, a similar argument with the isometry $\imath_2:(r,z,u) \mapsto (r, 2 a-z, u)$, allows us to regard $x_0$ as an interior point of $M_{t_0}$. 
\end{itemize}

If $g_0:=g(x_0,t_0) >1$, then  it has to satisfy $0 \le - \frac{k}{2} g_0^2  + (K_2 +K_3)\ g_0 $, so 
\be\lb{ubg}
g_0 \le \max\{\max_{M} g, 2 (K_2 +K_3)/k, 1\}.
\ee
Since $\varphi$  and $g$ are bounded independently of $T$, $ |L|^2$ is bounded on $[0,T[$ with a bound independent of $t$. 

Once we achieve the upper bound for  $|L|^2$, it  follows, like in  \cite{Hu84} and \cite{Hu86} (see \cite{Hu87}, \cite{Ath1} or \cite{CaMi1} for the volume-preserving version), that $|\nabla^jL|$ is bounded for every $j\ge 1$. If $T<\infty$, these bounds imply (cf. \cite{Hu84} pages 257, ff.) that   $X_t$ converges (as $t\to T$, in the $C^\infty$-topology) to a unique smooth  limit $X_T$. Now we can apply  the short time existence theorem to continue the solution after $T$,  contradicting the maximality of $[0,T[$.
\end{demo}

\section{Finer bounds for $\oH$}\label{bfoH}

 If $M$ is a revolution surface with boundary in the setting RSS2,  theorems \ref{graphbo} and \ref{t:preLTE}  tell us that  the flow $M_t$ of $M$ under \eqref{vpmf}, subject to the boundary conditions \eqref{bocon}, exists as long as $M_t$ does not touch the axis of revolution, and the generating curve of each $M_t$ is a graph over the same axis. This fact will allow us to improve the bounds for $\oH$ obtained in section \ref{broH} by using a Sturmian theorem (cf. \cite{Ang}), as it was done in \cite{Ath2} when $\oM$ is the Euclidean space. 

Before going on, we want to stress a couple of remarks which is necessary to have in mind from now on.

\begin{nota}\label{deriv_r}
After the results of sections \ref{broH} and \ref{graph}, it is clear that $r$ and $v$ remain bounded on $[0, t]$ for every $t\in[0,T[$. 
Then,  from \eqref{tNf} and the definition of $v$, it follows that $\dot{r}$ is also bounded. Now, from \eqref{k1f1} and \eqref{k2f}, we conclude that all the derivatives of $r$ are bounded if and only if $|\nabla^jL|$ is bounded for every $j\ge 0$, which was established within the proof of Theorem \ref{t:preLTE}.
\end{nota}

\begin{nota}\label{equiv_eq}
Since the solutions of equations \eqref{vpmf} and \eqref{vpmft} differ only on the parame\-trization of $M_t$, the functions  ${r}_t(z_t)$ and $r_t(z)$ corresponding to the solutions of  \eqref{vpmf} and \eqref{vpmft}, respectively, will have the same number of zeroes for their first derivatives (with respect to $z$) at the same values of $z$, and the same bounds for all their derivatives. Therefore, whenever we are interested on these magnitudes, we can use either equation \eqref{vpmf} or \eqref{vpmft} at our convenience. 
\end{nota}

\begin{lema}\label{discrt} Let $M_t$, with $t\in[0,T[$, be a maximal  solution of  \eqref{vpmf}  with initial condition $M$ in the setting RSS2 and satisfying \eqref{bocon}. For each $t\in[0,T[$, the set $Z_t=\{z\in[a,b]; \dot{r}_t =0\}$ is finite, the function $t \mapsto N(t):=\sharp(Z_t)$ is non increasing and, at the points $(z_0,t_0)$ satisfying $0  = \dot{r}_{t_0}(z_0) =\ddot{r}_{t_0}(z_0)$ there is a neighborhood where the number of zeroes decreases.
\end{lema}
\begin{demo} Thanks to Remark \ref{equiv_eq},  we can use equation \eqref{vpmft} for this proof.  In order to apply the Sturmian theorem (cf. Theorem 1.1 in \cite{Ath2}) to $\psi_t:=\dot{r}_t$, we compute $\parcial{\psi}{t}$ by taking the derivative of \eqref{tafl} with respect to $z$. The result is
\begin{align}
\parcial{\psi}{t}  = \frac{\ddot{\psi}}{\dot{r}^2 + f^2}  + &
\( \frac{\oH (f \ddot{r} -f' \dot{r}^2)}{ f^2 \sqrt{ \dot{r}^2+ f^2}}  -\( \frac{f'}{f} + (n-1)\frac{h'}{h}\)'   - 2 \frac{f' \ddot{r}}{f(\dot{r}^2 + f^2)} - 2 \ddot{r} \frac{ \ddot{r} +  f f' }{(\dot{r}^2 + f^2)^2} \right. \nn \\
& \left. - \frac{\dot{r}^2}{\dot{r}^2 + f^2} \(\frac{f'}{ f}\)' + 2 \dot{r}^2 \frac{f'}{f} \  \frac{\ddot{r} +  f f'}{(\dot{r}^2 + f^2)^2}\) \psi. \label{ddotr}
\end{align}

For every $t_1 < T$, let $\rho(t_1) = \min_{(x,t)\in M\times [0,t_1]} r_t(x) >0$. As follows from Remark \ref{deriv_r},  the coefficients of \eqref{ddotr} satisfy the hypotheses of the Sturmian theorem on $[0,t_1]$. Then, applying such theorem, we conclude that the statement is true for $t\in [0,t_1]$ and, since $t_1$ is arbitrary, also for all $t\in [0,T[$.
\end{demo}

A first consequence, deduced from the previous lemma, is that we can find an upper bound for the length of the generating curve of the evolving hypersurface. Indeed,
\begin{lema}\label{legc}
Let $M_t$, with $t\in[0,T[$, be a maximal  solution of  \eqref{vpmf}  with initial condition $M$ in the setting RSS2 and satisfying \eqref{bocon}.  There is a constant $c(N(0),r_2,\og, b - a)$ such that, for every $t\in[0,T[$, the length of the generating curve  of $M_t$ is lower than $c$.
\end{lema}
\begin{demo} 
From Lemma \ref{discrt}, $Z_t=\{a=z_1, ..., b=z_{N(t)}\}$ and $N(t)\le N(0)$. The computation of the length of the generating curve $\gamma_t$ gives 
$$L(\gamma_t)= \sum_{j=1}^{N(t)-1} \int_{z_j}^{z_{j+1}} \sqrt{\dot{r}(z)^2 + f(r(z))^2} \ dz \le f(r_2) (b-a) + (N(0)-1) r_2 =:c.$$
\end{demo}

Now we are in position to use the previous lemmas in order to achieve finer bounds for $\oH$. With \lq\lq finer" we mean that, unlike Corollary \ref{c:uboH}, here we obtain a double side bound for $\oH$ independent of the lower bound for $r$. In fact,
\begin{prop}\label{t:fboH}
Let $M_t$, with $t\in[0,T[$, be a maximal  solution of  \eqref{vpmf}  with initial condition $M$ in the setting RSS2 and satisfying \eqref{bocon}. Then we can find constants $h_i(V, \og,n,b-a,\vle(M),\mathcal{A}_V)> 0$, $i=1,2$   such that $h_1\le \oH \leq h_2$, where $\mathcal{A}_V$ is the infimum of the $n$-volume of the hypersurfaces satisfying the constraint of the setting RSS2 and enclosing a volume $V$.
\end{prop}
\begin{demo}
As in the proof of Proposition \ref{p:uboH} we have $I_1\ge 0$ and, starting with formula \eqref{I1pos2} for $I_1$, we can write:
\begin{align}
I_1 &\le \sigma\ \frac{(h^{n-1})'(r_2)}{\vle(M_t)}\ \frac{\pi}{2} \ L(\gamma_t) \le \sigma\ \frac{(h^{n-1})'(r_2)}{\mathcal{A}_V} \frac{\pi}{2} \ c \label{fbI2},
\end{align}
where the last inequality follows by Lemma \ref{legc}. On the other hand, from the expression  \eqref{I2pos} for $I_2$, we obtain
\begin{align}
I_2 \le \frac\sigma{\mathcal{A}_V} (b-a)\ \max_{r\in[0,r_2]}((n-1) h' f + f' h) h^{n-2})(r). \label{ubI_2}
\end{align}

To prove that $I_2$ has a positive uniform lower bound on $[0,T[$, let us consider 
\begin{equation} \label{defr3}
r_3 = \beta^{-1}\(\fracc{V}{2 (b-a) \sigma}\).
\end{equation}
Observe that, by the above definition, $0 < r_3 < r_1 < r_2$. 
 Let $L_3$ be the measure of the set $\mathcal{L}_3 = \{ z\in[a,b]; r_t(z)\ge r_3\}$. From \eqref{defr3}, \eqref{defr1} and \eqref{fVrz}, we have
 \begin{align}\frac{V}{2\sigma} &= \int_a^b  \int_{0}^{r_t(z)} f \ h^{n-1} \ dr  \ dz - \int_a^b  \int_0^{r_3}  f \ h^{n-1} \ dr \ dz = \int_a^b   \int_{r_3}^{r_t(z)} f \ h^{n-1} \ dr \ dz \nn\\ &
  \le \int_{\mathcal{L}_3}  \int_{r_3}^{r_t(z)} f \ h^{n-1} \ dr  \ dz \le f \ h^{n-1} (r_2) \ L_3\ (r_2-r_3), \lb{in1}
 \end{align}
which implies that $L_3\ge \rho_3>0$ for some $\rho_3$ independent of $t$. Now, we can use again \eqref{I2pos} to get
 $$I_2 \ge \frac{\sigma}{\vle(M_t)}  \int_{\mathcal{L}_3} \((n-1) h' f + f' h\) h^{n-2} dz \ge \frac{\sigma  \rho_3}{\vle(M)}\ \min_{r\in[r_3,r_2]}((n-1) h' f + f' h) h^{n-2})(r) .$$
 This gives the desired positive lower bound for $I_2$, which, together with the other inequalities obtained for $I_1$ and $I_2$, finishes the proof of the proposition.
 \end{demo}

\section{Discreteness of the first singularity set} \label{dis}

According to Theorem \ref{preLTE}, at the points $z$ where one has a singularity at time $T$ one must have $\lim_{t\to T}r_t(z)=0$. It seems that these points must be  limits of   zeroes of $\dot{r}_t$. We shall show that this is the case, which, thanks to Lemma \ref{discrt}, gives the discreteness of the set of singular points at time $T$.

By Lemma \ref{discrt}, $N(t)$ is finite and does not increase with time. Moreover, if $0  = \dot{r}_{t_0}(z_0) =\ddot{r}_{t_0}(z_0)$, there is a neighborhood where the number of zeroes decreases. Then, we can find some time $t_1\in [0,T[$ such that  for every $t\in[t_1,T[$, $N(t)$ is constant and $0   \ne\ddot{r}_{t}(z)$ for every $z\in Z_t$. Thus all the critical points of $r_t$  are local maxima or minima. Hereafter $\xi_i(t)$ denotes the minima and $\eta_i(t)$, the maxima.

Notice that, if the number of critical points for $r_t$ on $[a,b]$ is $2k+2$, we may have the following situations
\begin{align*}
a < \xi_1(t) < \eta_1(t) < ... < \xi_k(t) <\eta_k(t) < b\\
a  < \eta_1(t) < \xi_1(t)< ... <\eta_k(t) < \xi_k(t) < b
\end{align*}
and, if the number of critical points for $r_t$ is $2k+3$,
\begin{align*}
a < \xi_1(t) < \eta_1(t) < ... < \xi_k(t) <\eta_k(t) < \xi_{k+1}(t)< b\\
a  < \eta_1(t) < \xi_1(t)< ... <\eta_k(t) < \xi_k(t) < \eta_{k+1}(t)< b.
\end{align*}

On the other hand, as $\ddot{r}_t(\xi_i(t)) \ne 0 \ne \ddot{r}_t(\eta_i(t))$, by the implicit function theorem, the solutions $\xi_i(t)$, $\eta_i(t)$ of the equation $\dot{r}_t(z) = 0$ are smooth as functions of time. In addition, since the map $(r,z,u) \mapsto (r, 2 z_0- z,u)$ is an isometry of $\oM$ for every $z_0\in \re$, the same arguments given in Lemma 2.3 of \cite{Ath2} (together with those of 5.1 in \cite{AAG}) imply the existence of the limits 
$$\xi_i(T)=\lim_{t\to T} \xi_i(t) \qquad \text{ and } \qquad \eta_i(T)=\lim_{t\to T} \eta_i(t).$$ 

We can now state the main theorem of the present section.

\begin{teor} \lb{cdps} Let $M_t$, with $t\in[0,T[$, be a maximal  solution of  \eqref{vpmf}  with initial condition $M$ in the setting RSS2 and satisfying \eqref{bocon}.   For any closed interval $[c,d]$ not containing any of the points $a, \xi_i(T), \eta_i(T),b$, there exist $\delta>0$ and $t_2\in[0,T[$ such that $r_t(z)\ge \delta$ for $z\in [c,d]$, $t\in]t_2,T[$.
\end{teor}
\begin{demo}
First, assume that $[c,d]\subset]\xi_j(T),\eta_j(T)[ $ (the other cases can be proved in a similar way,  in a half of them with the inequalities reversed and using upper right Dini derivative instead of the lower one). Let $a',b'$ be real numbers satisfying $[c,d]\subset ]a',b'[$ and  $[a',b']\subset ]\xi_j(T),\eta_j(T)[$ .
By the definition of $t_1$ (at the beginning of this section), $\dot{r}_t(z)>0$ for $\xi_j(t)< z < \eta_j(t)$ and $t\in[t_1,T[$. Moreover, by the continuity of  $\xi_j(t)$ and $\eta_j(t)$ on $[t_1,T]$, there is  a $t_2\in]t_1,T[$ such that  $\dot{r}_t(z)>0$ for $(z,t)\in [a',b'] \times ]t_2,T[$.
On this domain, by discarding the positive addends in \eqref{ddotr}, we have
\begin{align}
\parcial{\psi}{t}  \ge \frac{\ddot{\psi}}{\dot{r}^2 + f^2}  \ +\  &
T(\psi)\ \psi  -\frac{f' \, \oH}{f^2 \sqrt{ \dot{r}^2+ f^2}} \psi^3,\lb{pomt}
\end{align}
being
\begin{align*}T(\psi) &= \frac{\oH \ddot{r} }{ f \sqrt{ \dot{r}^2+ f^2}}  -\frac{ f'' }{f} - 2 \frac{f' \  \ddot{r}}{f(\dot{r}^2 + f^2)} - 2 \frac{ \ddot{r}^2 +  f f' \ddot {r}}{(\dot{r}^2 + f^2)^2}  \ \nn \\
& - (n-1) \frac{h''  }{h} + \frac{\dot{r}^2}{f\(\dot{r}^2 + f^2\)} \( -f'' + \frac{2 \ddot{r}  f'}{\dot{r}^2 + f^2}\).
\end{align*}

Next, taking $A$ and $M_1$ two constants to be determined later, we define
\be\lb{wAM}
w(z,t) = e^{- A t} \sin(\mu (z-a')), \qquad \text{ where } \mu = \frac{\pi}{b'-a'} \quad \text{and }\quad A\ge\mu^2 + M_1.
\ee
Observe that $w$ satisfies
\be
\parcial{w}{t} \le - \(\mu^2  + M_1\) w \le - \frac{\mu^2}{f^2 + \dot{r}^2} w - M_1 w =  \frac{\ddot{w}}{f^2 + \dot{r}^2}  - M_1 w.\lb{pwt}
\ee

Let  us consider the difference $\frak{w}=\psi - w$. From \eqref{pomt} and \eqref{pwt}, we obtain the estimate:
\begin{equation}\label{pdif}
 \parcial{\frak{w}}{t}  \ge \frac{ \ddot{\frak{w}}}{\dot{r}^2 + f^2}  + 
\(T(\psi)-\frac{f' \, \oH}{f^2 \sqrt{ \dot{r}^2+ f^2}}\) \psi  + M_1 w +\frac{f' \, \oH}{f^2 \sqrt{ \dot{r}^2+ f^2}} (\psi-\psi^3).
\end{equation}

 Let us define the function $\frak{w}_m(t):= \min_{z\in[a',b']} \frak{w}(z,t)$. For any $t'\in[t_2,T[$, let $I_{t'} = \{z\in[a',b']; \frak{w}_m(t') = \frak{w}(z,t') \}$. 
 We have two possibilities for any fixed $\overline t \in [t_2,T[$:  
either $\{a',b'\} \cap I_{\overline t} \ne \emptyset$ or $\{a',b'\} \cap I_ 
{\overline t} = \emptyset$.

\begin{equation}\label{a'It}
\hspace*{-4.2cm} \text{If $a'$ or $b'$ are in $I_{\overline t}$, then $\begin{array}{l} \text{$\frak{w}(z, 
\overline t) \ge \frak{w}_m(\overline t) = \psi_{\overline t}(a') 
 >0$ or } \medskip\\
 \text{$ \frak{w}(z,\overline t) \ge \frak{\frak{w}}_m(\overline  
t) = \psi_{\overline t}(b')>0$.}\end{array}$}
\end{equation}
  
In the other case (that is, if $\{a',b '\} \cap I_{\bar t} = \emptyset$), let us call $]t_3, t_4[$ the maximal open interval in $[t_2,T]$ such that $\overline t \in ]t_3,t_4[$ and $\{a', b'\}\cap I_t = \emptyset$ for every $t\in]t_3,t_4[$ (notice that $t_3 = t_2$ or $\{a', b'\}\cap I_{t_3} \ne \emptyset$ and that $t_4=T$ or $\{a', b'\}\cap I_{t_4} \ne \emptyset$).

On this interval  we can write
\be\lb{pdifm}
\parcial{\frak{w}_m}{t} (t) \ge \min_{z\in I_t} 
\(\(T(\psi)-\frac{f' \, \oH}{f^2 \sqrt{ \dot{r}^2+ f^2}}\) \psi  + M_1 w \)+ \min_{z\in I_t}  \frac{f' \, \oH}{f^2 \sqrt{ \dot{r}^2+ f^2}} (\psi-\psi^3),
\ee
 where $\parcial{\frak{w}_m}{t} (t)$ must be understood in the sense of a lower right Dini derivative (cf. page 160 of  \cite{SoBoCh}).

Using the function $\Psi(t)=  \min_{z\in I_t} (\psi -\psi^3)(z,t)$ (defined on  $[t_3,t_4[$), we decompose 
\be
 [t_3,t_4[= J_1 \cup J_2, \text { where } J_1 = \Psi^{-1}(]0,\infty[ ) \text{ and } J_2 =  \Psi^{-1}(]-\infty, 0]). \nonumber
\ee 

For $t\in J_2$, there is a $z_t\in I_t$ such that $\psi(z_t,t) \ge 1$. Hence, for every $z\in[a',b']$, 
\be \lb{J2}
\frak{w}(z,t) \ge \frak{w}_m(t) = \psi(z_t,t) - w(z_t,t)\ge 1- w(z_t,t) > 0.
\ee 

If $t_0\in J_1$, let $(\tau,\eta[$ be the maximal interval containing $t_0$ and contained in $J_1$ (this interval is $]\tau,\eta[$ if $\tau > t_3$ and  it may be $[t_3,\eta[$ when $\tau= t_3$). Notice that  $\Psi(t)>0$  on $]\tau,\eta[$ and the function $t \mapsto  \min_{z\in I_t}  \frac{f' \, \oH}{f^2 \sqrt{ \dot{r}^2+ f^2}} (\psi-\psi^3)(z,t)$ is positive at the same $t$ as $\Psi(t)$. Therefore, by \eqref{pdifm}, we have
\be\lb{pdiftmi}
\parcial{\frak{w}_m}{t} (t) \ge \min_{z\in I_t} 
\(\(T(\psi)-\frac{f' \, \oH}{f^2 \sqrt{ \dot{r}^2+ f^2}}\)\ \psi  + M_1 w \)(z,t) \qquad \text{on} \quad ]\tau,\eta[.
\ee

On the other hand, at  $z\in I_t$, 
\be
0= \dot{\frak{w} }= \dot{\psi} - \dot{w} , \quad \text{ that is, }  \quad \dot{\psi} = \dot{w} = e^{-A t} \mu \cos(\mu (z-a')). \nonumber
\ee
Thus $|\dot{\psi}| =| \ddot{r}| \le \mu$ on $I_t$, which can be  used to bound
\begin{align}
&T(\psi) -\frac{f' \, \oH}{f^2 \sqrt{ \dot{r}^2+ f^2}} \nn \\
&\ge - \sup_{r\in[0,r_2]}  \( h_2 (\mu+f') +2 |f'' |+4 \mu  f'+ 2 \mu^2 + 2 f f' \mu + (n-1) \frac{h''  }{h}  \) = - C_5, \label{bTom}
\end{align}
where $C_5<\infty$  because $- \frac{h''}{h}$ is a sectional curvature of the ambient space (cf. \eqref{scoM}) and $f$ has all its derivatives bounded on a closed interval. If $\bar z \in I_t$ is the point where the minimum in the right hand side of \eqref{pdiftmi} is attained, then taking $M_1=C_5$ and substituting \eqref{bTom} in \eqref{pdiftmi},   we get
\be\label{lpfwt}
\parcial{\frak{w}_m}{t}(t) \ge - C_5 \  (\psi -w)(\bar z, t) = -C_5\   \frak{w}_m(t) \quad \text{ on } \quad ]\tau,\eta[.
\ee

From \eqref{lpfwt}  we conclude that, on $]\tau,\eta[$, 
\begin{align}
\frak{w}_m (t) &\ge e^{-C_5 (t-\tau)} \frak{w}_{m}(\tau) \ge 
e^{-C_5 (T-t_2)} \frak{w}_{m}(\tau). \lb{lbfw}
\end{align}
For every  $z\in I_\tau$,
\be\label{fraw-}
\frak{w}_m(\tau) = \psi(z,\tau) - e^{-A\tau} \sin(\mu(z-a')) \ge  \psi(z,\tau) - e^{-A t_2}.
\ee
We have two possibilities for $\tau$. Either $\tau=t_3$   or $\Psi(\tau)= 0$. 

In the first case ($\tau=t_3$),  as remarked after \eqref{a'It}, we have again two possibilities:
\begin{itemize}
\item $\{a',b'\} \cap I_{t_3} \ne \emptyset$, then $\frak w_m(\tau) =\ \frak w_m(t_3) = \left\{\begin{matrix} &\psi(a',t_3) >0\\ &\text{or} \\& \psi(b',t_3) >0\end{matrix}\right.$.
\item $t_3= t_2$, then $\Psi(t_2)\ge 0$ and, for $z\in I_{t_2}$, 
$$m_1:= \min_{[a',b']} \psi (\,\cdot \,,t_2) \leq \psi (z, t_2) \leq  1.$$
\end{itemize} 
Then, from \eqref{fraw-}, $\frak{w}_{m}(\tau)>0$ if $A \ge - \frac1{t_2} \ln m_1$, where  $m_1=\min_{z\in[a',b']}\dot{r}(z, t_2)$. 

In the second case, there is a $z_\tau \in I_\tau$ such that $1= \psi (z_\tau,\tau) $, thus $\frak{w}_{m}(\tau)>0$ if $A >0$.  

Since $A$ has also to satisfy \eqref{wAM}, we define $A = \max\{ - \frac1{t_2} \ln m_1, C_5+\mu^2\}$. With this election,  having also in mind the inequalities  \eqref{a'It},  \eqref{J2} and \eqref{lbfw}, we obtain $\frak{w}\ge 0$ on $[a',b']\times [t_2,T[$, i.e.
$\dot{r}(z,t) \ge e^{-At} \sin(\mu (z-a'))$ and, integrating with respect to $z$ between $a'$ and $z$, for $z\in[c,d]$,
\begin{align*}
r(z,t) & \ge  \frac{e^{-At}}{\mu} \(1- \cos\(\frac{\pi(z-a')}{b'-a'}\)\) 
\ge  \frac{e^{-AT}}{\mu} \(1- \cos\(\frac{\pi(c-a')}{b'-a'}\)\) = \delta >0,
\end{align*}
which finishes the proof of the theorem.
 \end{demo}

\section{Small ${ \mathbf n}$-volume for $\mathbf{ M}$ implies  long time existence and convergence.} \label{volp}

\begin{lema}\lb{smallv}In the setting RSS, if $\vle(M) \le \fracc{V}{(b-a)}\frac{\delta(r_1)}{\beta(r_1)}$, then there is a $r_0>0$ such that $r_t>r_0$ for every $t$ such that the flow $M_t$ exists.
\end{lema}
\begin{demo} If there is some $t_0$ and $x_0$ for which $r_{t_0}(x_0)=0$, using  \eqref{vle>+} and \eqref{fVr1}, we obtain
\be
\vle(M) \ge \vle(M_{t_0})> \sigma \int_0^{r_1}   h(r)^{n-1} dr  \ge \frac{V \delta(r_1)}{(b-a)\beta(r_1)},
\ee
in contradiction with the hypothesis.\end{demo}

\begin{nota}
Notice that $\ds \frac{\delta(r_1)}{\beta(r_1)} \leq 1$ with equality if $f(r) \equiv 1$. So the hypothesis on the $n$-volume of $M$  when $\oM$ is the Euclidean space coincides with that of \cite{Ath1}. For $n=2$ and $\oM$ a space of constant curcvature $\lambda$, the explicit value of the upper bound in Lemma \ref{smallv} is $\ds \frac{2 \pi}{-\lambda} \(-1 + \sqrt{1- \frac{\lambda V}{\pi (b-a)}}\)$.
\end{nota}

\begin{teor}\lb{fth}
If $M$ belongs to the setting RSS2 and has $\vle(M) \le \fracc{V}{(b-a)}\frac{\delta(r_1)}{\beta(r_1)}$, then the solution of \eqref{vpmf} satisfying \eqref{bocon}  is defined for all $t>0$ and converges to a revolution hypersurface of constant mean curvature in $\oM$.
\end{teor}
\begin{demo}
By Lemma \ref{smallv} and Theorem \ref{t:preLTE}, we know  that the solution of \eqref{vpmf} is defined on $[0,\infty[$. Using again Lemma \ref{smallv}, together with Proposition \ref{bor}, we conclude that $r$ is double side bounded uniformly on $[0,\infty[$. After the results of section \ref{graph}, it is clear that  $\dot{r}$ remains bounded all the time. 
In addition, the proof of Theorem \ref{t:preLTE} shows that  $|\nabla^jL|^2$ is uniformly bounded for every $j\ge 0$.  Once we have all these bounds, it follows from \eqref{k1f1} and \eqref{k2f} that all the derivatives of $r$ are bounded on $[0,\infty[$.  

We are now in position to  apply Arzel\`a-Ascoli Theorem to assure the existence of a sequence of maps  $ r_{t_i}$ satisfying \eqref{tafl} which $C^\infty$-converges to  a smooth  map  $r_\infty : [a,b]\flecha \re^+$ satisfying \eqref{tafl}. A standard argument like in \cite{CaMi1} proves that $ M_\infty =(z, r_\infty(z),u)$
has constant mean curvature.

On the other hand, for any $0<\tau_1 < \tau_2$, from \eqref{tafl}, we have 
\begin{align}
\max_M | r(x,\tau_2) - r(x,\tau_1)| & \le \max_M \int_{\tau_1}^{\tau_2} \left|\parcial{r}{t}(x,t)\right| dt \nn \\
&\le \max_M \int_{\tau_1}^{\tau_2}  |\oH - H| \sqrt{\dot{r}^2 + f^2} dt \le C \varepsilon(t) \ (\tau_2-\tau_1), \nn
\end{align}
where $\varepsilon(t) \to 0$ as $t\to \infty$. Then, arguing as in \cite{HuB} (end of section 2), we conclude that $r_t$ $C^\infty$-converges to   $r_\infty$.
Therefore, $ M_\infty =(z, r_\infty(z),u)$ is, up to tangential diffeomorphisms, the limit of $M_t$. 
\end{demo}

{\footnotesize

\bibliographystyle{alpha}

}

\vskip1truecm

{\small 

\begin{tabular}{ c c }
Address & Universidad de Valencia.\\
      \       & 46100-Burjassot (Valencia) Spain \\
      \       & email: esther.cabezas@uv.es\\
      \       & and\\
      \       & email: miquel@uv.es 
\end{tabular}

}

\end{document}